\documentclass[global]{svjour3}
\usepackage{amssymb}
\usepackage{amsfonts}
\usepackage{amsmath}
\usepackage{color}

\input{tcilatex}
\begin{document}

\journalname{Vietnamese Journal of Mathematics}
\title{\textbf{Robust generalized S-Procedure}}
\author{N. Dinh  \and M.A. Goberna \and D.H. Long  \and {\mbox{M. Volle}}}

\institute{
Nguyen Dinh,  Corresponding author  \at
	International University, Vietnam National University- HCMC,
	Linh Trung ward, Thu Duc district, Ho Chi Minh city, Vietnam 
           ({\it ndinh@hcmiu.edu.vn}); 
	Vietnam National University - HCMC,
 	Linh Trung ward, Thu Duc district, Ho Chi Minh city, Vietnam.
\and M.A. Goberna \at
	Department of Mathematics, University of Alicante, 03080 Alicante, Spain 
	 (\textit{mgoberna@ua.es}).
\and D.H. Long \at
	Ho Chi Minh City University of Technology, Vietnam National University - HCMC, 
	268 Ly Thuong Kiet, Ward 14, District 10, Ho Chi Minh city, Vietnam
	({\it danghailong@hcmut.edu.vn}).
\and M. Volle  \at
            D\'{e}partement de Math\'{e}matiques, LMA EA 2151, Avignon, France 
	({\it michel.volle@univ-avignon.fr}).
}

\maketitle

\noindent{\sl Dedicated to Professor Hoang Xuan Phu on the occasion of his  $70^{th}$ birthday.} 

\bigskip

\begin{abstract}
We introduce in this paper the so-called robust generalized $S$-procedure
associated with a given robust optimization problem. We provide a primal
characterization for the validity of this procedure as well as a dual
characterization under the assumption that the decision space is locally
convex. We also analyze an extension of the mentioned robust $S$-procedure
that incorporates a\textit{\ }right-hand side function.
\end{abstract}


%


\renewcommand{\proofname}{Proof}


{\subclass{90C17 \and 90C46 \and 49N15 \and 46N10}

\section{Introduction}

\label{sec:1}

According to \cite{GL06},\ the ``$S$-procedure" is the name given in \cite%
{AG63} - a monograph on the stability of nonlinear automatic control
systems\ where certain function denoted by $S$ (the initial of ``stability")
plays a crucial role - to a method which allows to solve systems of
quadratic inequalities by solving suitable linear matrix inequality
relaxations. Yakubovich was the first one to give sufficient conditions for
these relaxations to be exact in the sense that it is possible to obtain a
solution for these quadratic systems using their corresponding relaxations 
\cite{Yabu71}. The classic theorem of losslessness of the $S$-procedure (in
short $S$-procedure) characterizes when certain concave quadratic inequality
is consequence of a given system of convex quadratic inequalities in terms
of the concave function being a nonnegative linear combination of the convex
constraints and a suitable trivial constraint. This result allows to get
optimality conditions for convex quadratic problems. In fact, given a family
of convex quadratic functions $q,q_{1},...,q_{m},$ and 
\begin{equation*}
\overline{x}\in X:=\left\{ x\in \mathbb{R}^{n}:q_{i}\left( x\right) \leq
0,i=1,...,m\right\} ,
\end{equation*}
{observing that} $\overline{x}$ is an optimal solution of $\inf\limits_{x\in X}q\left(
x\right) $ if and only if the concave inequality $q\left( \overline{x}%
\right) -q\left( x\right) \leq 0$ is consequence of the system of convex
quadratic inequalities $\left\{ q_{i}\left( x\right) \leq
0,\; i=1,\ldots,m\right\} .$ The classic $S$-procedure can be seen as a nonlinear
version of Farkas' lemma and has a clear algebraic nature in the sense that
it is closely related with convex analysis and the Lagrange and Fenchel
duality theories for optimization problems.

It is worth observing that certain authors call $S$-lemma to the theorem of
the $S$-procedure when there exists a unique constraint \cite{DP06}, while
other authors prefer to call $S$-lemma to the characterization of the
concave quadratic inequalities which are consequence of systems of convex
quadratic inequalities and $S$-procedure to its application to nonlinear
automatic control systems \cite{PT07}.

The validity of the $\QTR{frametitle}{S}_{F}$\textit{-}procedure provides
optimality conditions for deterministic optimization problems, but it is
well-known that the data are uncertain in many optimization problems arising
in practice. For instance, in robust optimization the conservative
decision-maker assumes the existence of an uncertainty set $U$ and a family
of uncertain functions $f,f_{1},...,f_{m}:\mathbb{R}^{n}\times
U\longrightarrow \mathbb{R},$ and considers elements of the robust feasible
set 
\begin{equation*}
X:=\left\{ x\in \mathbb{R}^{n}:f_{i}\left( x,u\right) \leq 0,i=1,...,m,u\in
U\right\}
\end{equation*}%
in order to guarantee the feasibility of the chosen decision under any
conceivable scenario. Then, the (pessimistic) robust optimization problem to
be solved reads 
\begin{equation*}
\inf\limits_{x\in X}\sup\limits_{u\in U}f\left( x,u\right) .
\end{equation*}

Let us show how uncertain quadratic systems arise in astrophysics. Descartes 
\cite{Descartes1644} conceived in 1644 the tesselation of the space into
regions of \emph{influence} of the fixed stars according to proximities for
the Euclidean distance. In his model, if $T$ is the (huge) set of fixed
stars, the region of influence of certain $s\in T$ is 
\begin{equation*}
\begin{array}{ll}
\left\{ x\in \mathbb{R}^{3}:\left\Vert x-s\right\Vert \leq \left\Vert
x-t\right\Vert ,\text{ }t\in T\right\} & =\left\{ x\in \mathbb{R}%
^{3}:\left\Vert x-s\right\Vert ^{2}\leq \left\Vert x-t\right\Vert ^{2},\text{
}t\in T\right\} \\ 
& =\left\{ x\in \mathbb{R}^{3}:\left\langle \left( t-s\right)
,x\right\rangle \leq \tfrac{\left\Vert t\right\Vert ^{2}-\left\Vert
s\right\Vert ^{2}}{2},t\in T\right\} ,%
\end{array}%
\end{equation*}%
a polyhedron presently called Voronoi cell of $s$ (see a 2D representation
of Descartes' tesselation\ in \cite[Figure 2]{LP12}). This model became
obsolete after the publication of Newton's Principia in 1687 \cite%
{Newton1687}, where he stated the law of universal gravitation that forces
to replace the distance $\left\Vert x-t\right\Vert $ in the above cell by
the gravitational potential created by the star $t$ at $x,$ i.e., $%
u_{t}\left\Vert x-t\right\Vert ^{-2},$ for $x\neq t,$\ where $u_{t}$ is the
product of the gravitational constant times the mass of star $t.$ Assume
that this uncertain scalar ranges on some interval $U_{t}$ and define the
uncertainty set $U:=\prod\limits_{t\in T}U_{t}.$ For $x\notin T$ one has%
\begin{equation*}
u_{s}\left\Vert x-s\right\Vert ^{-2}\geq u_{t}\left\Vert x-t\right\Vert
^{-2}\Longleftrightarrow u_{t}\left\Vert x-t\right\Vert ^{2}\leq
u_{s}\left\Vert x-s\right\Vert ^{2},
\end{equation*}%
so that the (uncertain) influence region of $s\in T$ writes 
\begin{equation*}
\left\{ x\in \mathbb{R}^{3}:q_{t}\left( x,u\right) \leq 0,\text{ }t\in
T\right\} ,
\end{equation*}
where $q_{t}\left( x,u\right) =u_{t}\left\Vert x-t\right\Vert
^{2}-u_{s}\left\Vert x-s\right\Vert ^{2}$ is a non-convex (actually D.C.)
quadratic function for all $t\in T.$

The deterministic classic $S$-procedure has been extended in different ways
to more general functions, see for instance the survey papers \cite{DP06}, 
\cite{GL06} and \cite{PT07}. This can be done, e.g., in either a geometric
way (\cite{Lara&JBHU}, where the \textit{decision set} $X$ is the solution
set of a given constraint system posed in $\mathbb{R}^{n}$) or in the
following perturbational way (\cite{Lara-23}, \cite{Rock74}, \cite{RW-21}, 
\cite{VB-24}, etc., where the decision set $X$ is arbitrary): Let $Y$ be a
locally convex Hausdorff topological vector space (lcHtvs in short) of
perturbations whose null element is $0_{Y}$ and whose topological dual is $%
Y^{\ast },$ with duality product $\langle \lambda ,y\rangle :=\lambda (y)$
for all $(\lambda ,y)\in Y^{\ast }\times Y.$

Let 
\begin{equation*}
F\colon X\times Y\longrightarrow \overline{\mathbb{R}}:=\mathbb{R\cup }%
\left\{ \pm \infty \right\}
\end{equation*}%
be a perturbation functional, also called \textit{Rockafellian} in \cite%
{Lara-23} and \cite{RW-21}, and let 
\begin{equation*}
\inf_{x\in X}F(x,0_{Y})
\end{equation*}%
be the corresponding optimization problem. The \textit{generalized }$S$%
-procedure associated with $F$\ introduced in \cite{VB-24}, $%
\QTR{frametitle}{S}_{F}$\textit{-procedure} in brief, is concerned with the
next two statements: 
\begin{align*}
& (\mathbb{A})\quad F(x,0_{Y})\geq 0\text{ for all }x\in X, \\
& (\mathbb{B})\quad \exists \bar{\lambda}\in Y^{\ast }:F(x,y)+\langle \bar{%
\lambda},y\rangle \geq 0\text{ for all }(x,y)\in X\times Y.
\end{align*}%
Then, the $\QTR{frametitle}{S}_{F}$\textit{-}procedure is said to be\textit{%
\ valid} if $(\mathbb{A})\QTR{frametitle}{\Longleftrightarrow }\left( 
\mathbb{B}\right) $ or, equivalently, $\left( \mathbb{A}\right)
\Longrightarrow (\mathbb{B})$.

Roughly speaking, we are concerned with the adaptation of the $%
\QTR{frametitle}{S}_{F}$\textit{-}procedure to robust optimization problems.
More precisely, in this paper we consider a given uncertainty set $U$ and a
corresponding family of Rockafellians $F_{u}\colon X\times Y\rightarrow 
\overline{\mathbb{R}},$ $u\in U.$ Our purpose is to introduce and
characterize a robust generalized\textit{\ }$S$-procedure associated with
the robust optimization problem 
\begin{equation*}
\inf_{x\in X}\left\{ \sup_{u\in U}F_{u}(x,0_{Y})\right\} ,
\end{equation*}%
called \textit{robust} $S_{\left\{ F_{u}\right\} _{u\in U}}$\textit{%
-procedure}, which is concerned with the statements 
\begin{align*}
& (\mathfrak{A})\quad \sup\limits_{u\in U}F_{u}(x,0_{Y})\geq 0\text{ for all 
}x\in X, \\
& (\mathfrak{B})\quad \exists (\bar{u},\bar{\lambda})\in U\times Y^{\ast
}:F_{\bar{u}}(x,y)+\langle \bar{\lambda},y\rangle \geq 0\text{ for all }%
(x,y)\in X\times Y.
\end{align*}%
Observe that, if $(\mathfrak{B})$ holds and $x\in X$, then 
\begin{equation*}
\sup\limits_{u\in U}F_{u}(x,0_{Y})\geq F_{\bar{u}}(x,0_{Y})\geq -\langle 
\bar{\lambda},0_{Y}\rangle =0,
\end{equation*}%
i.e., $(\mathfrak{A})$ holds, too. Then, the robust $S_{\left\{
F_{u}\right\} _{u\in U}}$\textit{-}procedure\textit{\ }is said to be\textit{%
\ valid} if $(\mathfrak{A})\Longrightarrow (\mathfrak{B})$ or, equivalently, 
$(\mathfrak{A})\Longleftrightarrow (\mathfrak{B})$.

We also consider the \textit{robust generalized }$S$-procedure\textit{\ with
right-hand side (RHS in brief) function\ }$h\in \Gamma (X)$ (the family of
proper convex and lower semicontinuous functions from $X$ to $\overline{%
\mathbb{R}}$), \textit{robust} $S_{\left\{ F_{u}\right\} _{u\in U}}^{h}$%
\textit{-procedure} \textit{for short, }whose corresponding statements are 
\begin{align*}
& (\mathfrak{A}_{h})\quad \sup\limits_{u\in U}F_{u}(x,0_{Y})\geq h(x)\text{
for all }x\in X, \\
& (\mathfrak{B}_{h})\quad 
\begin{cases}
\forall a^{\prime }\in \limfunc{dom}h^{\ast },\text{ }\exists (\bar{u},\bar{y%
})\in U\times Y^{\ast }\text{ such that} \\ 
F_{\bar{u}}(x,y)+\langle \bar{\lambda},y\rangle \geq \langle a^{\prime
},x\rangle -h^{\ast }(a^{\prime }),\;\forall (x,y)\in X\times Y,%
\end{cases}%
\end{align*}%
where $h^{\ast }$ is the Fenchel conjugate of $h$ and $\limfunc{dom}h^{\ast
} $ is its domain. We say that the robust $S_{\left\{ F_{u}\right\} _{u\in
U}}^{h}$-procedure is\textit{\ valid} if $(\mathfrak{A}_{h})\QTR{frametitle}{%
\Longleftrightarrow }(\mathfrak{B}_{h})$ or, equivalently, $(\mathfrak{A}%
_{h})\Longrightarrow (\mathfrak{B}_{h})$ (the proof of $(\mathfrak{B}%
_{h})\Longrightarrow (\mathfrak{A}_{h})$ is postponed to Lemma \ref{lem:4-1}%
).

The rest of the paper is organized in three sections. We present in Section
2 a primal characterization for the validity of the robust $S_{\left\{
F_{u}\right\} _{u\in U}}$-procedure (Theorem \ref{thm:2-1} and Corollaries %
\ref{cor:2-1} and \ref{cor:2-2}). Section 3 analyzes the case when $X$ is a
lcHtvs (as $Y$), obtaining a dual characterization of the robust $S_{\left\{
F_{u}\right\} _{u\in U}}$-procedure (Theorem \ref{thm:3-1} and Corollary \ref%
{cor:3-1}). Finally, we characterize in Section 4 the validity of the robust 
$S_{\left\{ F_{u}\right\} _{u\in U}}^{h}$-procedure (Theorem \ref{thm:4-1}
and Corollary \ref{cor:4-1}).

We basically follow the notation of \cite{Z-02}, which is also the reference
for the results on convex analysis in lcHtvs, specially conjugate calculus.

\section{A primal characterization of the robust $S_{\left\{ F_{u}\right\}
_{u\in U}}$-procedure}

In this section $X$ is a given set and $\left\{ F_{u}\right) _{u\in U}$ is a
family of Rockafellians 
\begin{equation*}
F_{u}\colon X\times Y\longrightarrow \overline{\mathbb{R}}=\mathbb{R}\cup
\{+\infty \}\cup \{-\infty \},\quad \forall u\in U.
\end{equation*}

Given $G\colon X\times Y\longrightarrow \overline{\mathbb{R}}$ we denote by 
\begin{equation*}
\mathcal{G}:=\left\{ (y,r)\in Y\times \mathbb{R}:\exists x\in X,\;G(x,y)\leq
r\right\}
\end{equation*}%
the projection of the epigraph of $G$, 
\begin{equation*}
\limfunc{epi}G:=\left\{ (x,y,r)\in X\times Y\times \mathbb{R}:G(x,y)\leq
r\right\} ,
\end{equation*}%
on $Y\times \mathbb{R}$ and by 
\begin{equation*}
\limfunc{dom}G:=\left\{ (x,y)\in X\times Y:G(x,y)<+\infty \right\}
\end{equation*}%
the domain of $G,$ that is the projection of $\limfunc{epi}G$ on $X\times Y$.

We set $\mathbb{R}_{+}:=[0,+\infty \lbrack $ and $\mathbb{R}_{+}^{\ast
}=]0,+\infty \lbrack .$

\begin{lemma}
\label{lem:2-1} Given $G\colon X\times Y\longrightarrow \overline{\mathbb{R}}
$, the following statements are equivalent:

\textup{(i)} $G(x,0_Y)\ge 0,\; \forall x\in X$,

\textup{(ii)} $(0_Y,-1)\notin\mathbb{R}_+^\ast \mathcal{G}$,

\textup{(iii)} $(0_Y,-1)\notin\mathbb{R}_+ \mathcal{G}$.
\end{lemma}

\textbf{Proof} \textquotedblleft Non (i)\textquotedblright\ amounts to say
that%
\begin{equation*}
\begin{array}{ll}
\exists x\in X:G(x,0_{Y})<0 & \Longleftrightarrow \;\exists (x,\theta )\in
X\times \mathbb{R}_{+}^{\ast }:G(x,0_{Y})\leq -\theta \\ 
& \Longleftrightarrow \;\exists \theta \in \mathbb{R}_{+}^{\ast
}:(0_{Y},-\theta )\in \mathcal{G} \\ 
& \Longleftrightarrow \;\exists \theta \in \mathbb{R}_{+}^{\ast
}:(0_{Y},-1)\in \frac{1}{\theta }\mathcal{G} \\ 
& \Longleftrightarrow \;(0_{Y},-1)\in \mathbb{R}_{+}^{\ast }\mathcal{G},%
\end{array}%
\end{equation*}
that is \textquotedblleft Non (ii)\textquotedblright\ or, equivalently,
\textquotedblleft Non (iii)\textquotedblright .\hfill $\square \medskip $

For a subset $C$ of a lcHtvs we denote by $\limfunc{co}C$ its convex hull
and by $\overline{C}$ its closure.

\begin{lemma}
\label{lem:2-2} (see \cite[Proposition 2.1]{VB-24}) Let $G\colon X\times
Y\longrightarrow \overline{\mathbb{R}}$. The following statements are
equivalent

\textup{(i)} $\exists \lambda \in Y^{\ast }\text{ such that }G(x,y)+\langle
\lambda ,y\rangle \geq 0,\;\forall (x,y)\in X\times Y$,

\textup{(ii)} $(0_{Y},-1)\notin \overline{\limfunc{co}}\left( \mathbb{R}_{+}%
\mathcal{G}\right) $.
\end{lemma}

\textbf{Proof }Note that if $\limfunc{dom}G=\emptyset $ then (i) and (ii)
are both satisfied. So we can assume that $\limfunc{dom}G\neq \emptyset $,
which entails $\mathcal{G}\neq \emptyset .$

$\left[ \text{(i)}\Longrightarrow \text{(ii)}\right] $ By (i) there exists $%
\lambda \in Y^{\ast }$ such that 
\begin{equation*}
r+\langle \lambda ,y\rangle \geq 0\text{ for all }(y,r)\in \mathcal{G},
\end{equation*}%
that means $\mathcal{G}\subset \limfunc{epi}(-\lambda )$. Now $\limfunc{epi}%
(-\lambda )$ is a closed convex cone containing $\mathcal{G}$. By this we
have $\overline{\limfunc{co}}\left( \mathbb{R}_{+}\mathcal{G}\right) \subset 
\limfunc{epi}(-\lambda )$. Since $(0_{Y},-1)\notin \limfunc{epi}(-\lambda )$
it follows that (ii) holds.

$\left[ \text{(ii)}\Longrightarrow \text{(i)}\right] $ By the Hahn-Banach
separation Theorem (\cite[Theorem 3.4]{Rudin-91}), there exists $(\mu ,s)\in
Y^{\ast }\times \mathbb{R}$ such that 
\begin{equation*}
-s=\langle \mu ,0_{Y}\rangle -s>\alpha :=\sup_{\substack{ \theta \geq 0  \\ %
G(x,y)\leq r}}\theta \left( \langle \mu ,y\rangle +rs\right) .
\end{equation*}%
We have necessarily 
\begin{equation*}
G(x,y)\leq r\;\Longrightarrow \;\langle \mu ,y\rangle +rs\leq 0.
\end{equation*}%
Therefore $\alpha \leq 0$ and, since $\limfunc{dom}G\neq \emptyset $, $%
\alpha =0$. Thus $s<0$ and, setting $\lambda =\mu /s$, we get that 
\begin{equation*}
G(x,y)\leq r\;\Longrightarrow \;\langle \lambda ,y\rangle +r\geq 0
\end{equation*}%
and (i) holds.\hfill $\square \medskip $

\begin{theorem}
\label{thm:2-1} Let $F_{u}\colon X\times Y\longrightarrow \overline{\mathbb{R%
}}$, $u\in U$. The following statements are equivalent:

\textup{(i)} The robust $S_{\{F_{u}\}_{u\in U}}$-procedure is valid,\medskip

\textup{(ii)} $(0_{Y},-1)\notin \left[ \bigcap\limits_{u\in U}\overline{%
\limfunc{co}}\left( \mathbb{R}_{+}\mathcal{F}_{u}\right) \right] \setminus %
\left[ \mathbb{R}_{+}\bigcap\limits_{u\in U}\mathcal{F}_{u}\right] $.
\end{theorem}

\textbf{Proof }$\left[ \text{(i)}\Longrightarrow \text{(ii)}\right] $ Assume
that $(0_{Y},-1)\notin \mathbb{R}_{+}\bigcap\limits_{i\in U}\mathcal{F}_{u}$%
, where%
\begin{equation*}
\mathcal{F}_{u}:=\left\{ (y,r)\in Y\times \mathbb{R}:\exists x\in
X,\;F_{u}(x,y)\leq r\right\} .
\end{equation*}

We have to prove that $(0_{Y},-1)\notin \bigcap\limits_{u\in U}\overline{%
\limfunc{co}}\left( \mathbb{R}_{+}\mathcal{F}_{u}\right) $. Setting $%
G=\sup\limits_{u\in U}F_{u}$ we have that 
\begin{equation*}
\mathcal{G}=\bigcap_{u\in U}\mathcal{F}_{u},\text{ and }(0_{Y},-1)\notin 
\mathbb{R}_{+}\mathcal{G}.
\end{equation*}%
By Lemma \ref{lem:2-1} we obtain that 
\begin{equation*}
\sup_{u\in U}F_{u}(x,0_{Y})\geq 0\text{ for all }x\in X.
\end{equation*}%
Since the robust $S_{\{F_{u}\}_{u\in U}}$-procedure is valid, there exists $(%
\bar{u},\bar{\lambda})\in U\times Y^{\ast }$ such that 
\begin{equation*}
F_{\bar{u}}(x,y)+\langle \bar{\lambda},y\rangle \geq 0\text{ for all }%
(x,y)\in X\times Y.
\end{equation*}%
Applying Lemma \ref{lem:2-2} for $G=F_{\bar{u}}$ we obtain that 
\begin{equation*}
(0_{Y},-1)\notin \overline{\limfunc{co}}\left( \mathbb{R}_{+}\mathcal{F}_{%
\bar{u}}\right) .
\end{equation*}

$\left[ \text{(ii)}\Longrightarrow \text{(i)}\right] $ Assume that $(%
\mathfrak{A})$ holds. We have to check that $(\mathfrak{B})$ holds. Setting $%
G:=\sup\limits_{u\in U}F_{u}$, we have $\mathcal{G}=\bigcap_{u\in U}\mathcal{%
F}_{u}$ and 
\begin{equation*}
G(x,0_{Y})\geq 0\text{ for all }x\in X.
\end{equation*}%
By Lemma \ref{lem:2-1} we obtain that 
\begin{equation*}
(0_{Y},-1)\notin \mathbb{R}_{+}\bigcap\limits_{u\in U}\mathcal{F}_{u}.
\end{equation*}%
By (ii) it follows that 
\begin{equation*}
(0_{Y},-1)\notin \bigcap\limits_{u\in U}\overline{\limfunc{co}}\left( 
\mathbb{R}_{+}\mathcal{F}_{u}\right) .
\end{equation*}%
Then, by Lemma \ref{lem:2-2}, we obtain that $(\mathfrak{B})$ holds.\hfill $%
\square \medskip $

In order to recover well-known results on the $S_{F}$-procedure, we recall
that a subset $C$ of some lcHtvs space is said to be closed (resp. closed
convex) regarding another subset $D$ of the same space if $\overline{C}\cap
D=$\ $C\cap D$\ (resp., $\left( \overline{\limfunc{co}}C\right) \cap D=$\ $%
C\cap D$), see \cite{Bot-10} (resp., \cite{EV16}).

\begin{corollary}
\label{cor:2-1} (\cite[Theorem 2.1]{VB-24}) Let $F\colon X\times
Y\longrightarrow \overline{\mathbb{R}}$. The following statements are
equivalent:

\textup{(i)} The $S_{F}$-procedure is valid,

\textup{(ii)} $\mathbb{R}_{+}\mathcal{F}$ is closed convex regarding $%
\left\{ (0_{Y},-1)\right\} $.
\end{corollary}

\textbf{Proof }Let $U=\{u\}$ be a singleton and $F_{u}:=F$. The robust $%
S_{\{F_{u}\}_{u\in U}}$-procedure is valid if and only if the $S_{F}$%
-procedure is valid. By Theorem \ref{thm:2-1} this holds if and only if $%
(0_{Y},-1)\notin \left[ \overline{\limfunc{co}}\left( \mathbb{R}_{+}\mathcal{%
F}\right) \right] \setminus \left[ \mathbb{R}_{+}\ \mathcal{F}\right] $ that
is, if and only if, 
\begin{equation*}
\left[ \overline{\limfunc{co}}\left( \mathbb{R}_{+}\mathcal{F}\right) \right]
\cap \left\{ (0_{Y},-1)\right\} =\left[ \mathbb{R}_{+}\mathcal{F}\right]
\cap \left\{ (0_{Y},-1)\right\} ,
\end{equation*}%
namely $\mathbb{R}_{+}\mathcal{F}$ is closed convex regarding $\left\{
0_{Y},-1)\right\} $.\hfill $\square \medskip $

\begin{corollary}
\label{cor:2-2} (\cite[Corollary 2.1]{VB-24}) Let $F\colon X\times
Y\longrightarrow \overline{\mathbb{R}}$ be such that $\overline{\mathbb{R}%
_{+}\mathcal{F}}$ is convex. The following statements are equivalent:

\textup{(i)} The $S_{F}$-procedure is valid,

\textup{(ii)} $\mathbb{R}_{+}\mathcal{F}$ is closed regarding $\left\{
(0_{Y},-1)\right\} $.
\end{corollary}

\textbf{Proof }Since $\overline{\mathbb{R}_{+}\mathcal{F}}$ is convex we
have that 
\begin{equation*}
\overline{\limfunc{co}}\mathbb{R}_{+}\mathcal{F}=\overline{\mathbb{R}_{+}%
\mathcal{F}}.
\end{equation*}%
Therefore $\mathbb{R}_{+}\mathcal{F}$ is closed convex regarding $\left\{
(0_{Y},-1)\right\} $ if and only if $\mathbb{R}_{+}\mathcal{F}$ is closed
regarding $\left\{ (0_{Y},-1)\right\} $.\hfill $\square $


\section{Dual characterization of the robust generalized $S$-procedure}

From now on $X$ is a lcHtvs (as $Y$) with topological dual $X^{\ast }$,
whose null vectors are $0_{X}$ and $0_{X^{\ast }},$ respectively. The only
topology we consider in $X^{\ast }$ (resp. $X^{\ast }\times \mathbb{R}$) is
the weak*-topology.

Given $F_{u}\colon X\times Y\longrightarrow \overline{\mathbb{R}}$ for all $%
u\in U$, define 
\begin{align*}
p(x)& :=\sup_{u\in U}F_{u}(x,0_{Y}),\;\forall x\in X, \\
q(x^{\prime })& :=\inf_{(u,\mu )\in U\times Y^{\ast }}\left( F_{u}\right)
^{\ast }(x^{\prime },\mu ),\;\forall x^{\prime }\in X^{\ast },
\end{align*}%
and consider the projection of $\bigcup\limits_{u\in U}\limfunc{epi}%
(F_{u})^{\ast }$ on $X^{\ast }\times \mathbb{R}$, namely 
\begin{equation*}
\mathcal{F}^{\#}=\left\{ (x^{\prime },s)\in X^{\ast }\times \mathbb{R}%
:\exists (u,\mu )\in U\times Y^{\ast },\;(F_{u})^{\ast }(x^{\prime },\mu
)\leq s\right\} .
\end{equation*}%
One has (see, e.g.,  \cite{Lara&JBHU}, \cite{DGLV-19}, \cite
{DGLV-17},  \cite{DGV-around-24}, \cite{DGV-24}, \cite{DMVV-17}) 
\begin{gather*}
q^{\ast }(x)=\sup_{u\in U}(F_{u})^{\ast \ast }(x,0_{Y})\leq p(x)\text{ for
all }x\in X, \\
\overline{\limfunc{co}}\mathcal{F}^{\#}\subset \limfunc{epi}q^{\ast \ast
}\subset \limfunc{epi}p^{\ast }, \\
\limfunc{dom}q^{\ast }=\left\{ a\in X:\sup_{u\in U}(F_{u})^{\ast \ast
}(a,0_{Y})<+\infty \right\} , \\
\text{and, if }\limfunc{dom}q^{\ast }\neq \emptyset ,\;\limfunc{epi}q^{\ast
\ast }=\overline{\limfunc{co}}\mathcal{F}^{\#}.
\end{gather*}

Note that the statements $(\mathfrak{A})$, $(\mathfrak{B})$ read as follows: 
\begin{align*}
& (\mathfrak{A})\quad (0_{X^{\ast }},0)\in \limfunc{epi}p^{\ast }, \\
& (\mathfrak{B})\quad (0_{X^{\ast }},0)\in \mathcal{F}^{\#},
\end{align*}%
which shows that $(\mathfrak{B})\Longrightarrow (\mathfrak{A})$\ derives
from the inclusion $\mathcal{F}^{\#}\subset \limfunc{epi}p^{\ast }.$

Consider now the statements 
\begin{align*}
& (\mathfrak{A}^{\ast \ast })\quad \sup\limits_{u\in U}(F_{u})^{\ast \ast
}(x,0_{Y})\geq 0\text{ for all }x\in X, \\
& (\mathfrak{B}^{\ast \ast })\quad \exists (\bar{u},\bar{\lambda})\in
U\times Y^{\ast }:(F_{\bar{u}})^{\ast \ast }(x,y)+\langle \bar{\lambda}%
,y\rangle \geq 0\text{ for all }(x,y)\in X\times Y.
\end{align*}%
Note that $(\mathfrak{A}^{\ast \ast })$ reads 
\begin{equation*}
(\mathfrak{A}^{\ast \ast })\quad (0_{X^{\ast }},0)\in \limfunc{epi}q^{\ast
\ast }.
\end{equation*}%
Note that, since $\left( (F_{\bar{u}})^{\ast \ast }\right) ^{\ast }=(F_{\bar{%
u}})^{\ast }$, $(\mathfrak{B}^{\ast \ast })$ coincides with $(\mathfrak{B})$.

\begin{theorem}
\label{thm:3-1} Consider the following statements:

\textup{(i)} The robust $S_{\{F_{u}\}_{u\in U}}$-procedure is valid.

\textup{(ii)} The robust $S_{\{(F_{u})^{\ast \ast }\}_{u\in U}}$-procedure
is valid.

\textup{(iii)} $\mathcal{F}^\#$ is weak*-closed convex regarding $%
\{(0_{X^\ast},0)\}.$

We have $(\mathrm{i})\Longrightarrow (\mathrm{ii})\Longrightarrow (\mathrm{%
iii})$.\newline
If, moreover, 
\begin{equation*}
(H_{1})\qquad \qquad \inf_{x\in X}\sup_{u\in U}F_{u}(x,0_{Y})=\inf_{x\in
X}\sup_{u\in U}(F_{u})^{\ast \ast }(x,0_{Y})\neq +\infty
\end{equation*}%
holds, then $(\mathrm{i})\QTR{frametitle}{\Longleftrightarrow }(\mathrm{ii})%
\QTR{frametitle}{\Longleftrightarrow }(\mathrm{iii})$.
\end{theorem}

\textbf{Proof }$\left[ (\mathrm{i})\Longrightarrow (\mathrm{ii})\right] $
Assume (i).\ Since $(F_{u})^{\ast \ast }\leq F_{u}$ we have $(\mathfrak{A}%
^{\ast \ast })\Longrightarrow (\mathfrak{A})$. By this and the fact that $(%
\mathfrak{B}^{\ast \ast })\Longleftrightarrow (\mathfrak{B})$, we obtain
(ii).

$\left[ (\mathrm{ii})\Longrightarrow (\mathrm{iii})\right] $ Let $%
(0_{X^{\ast }},0)\in \overline{\limfunc{co}}\mathcal{F}^{\#}$. We have $%
(0_{X^{\ast }},0)\in \limfunc{epi}q^{\ast \ast }$. This is $(\mathfrak{A}%
^{\ast \ast })$. By (ii), $(\mathfrak{B}^{\ast \ast })$, alias $(\mathfrak{B}%
)$, holds, that is $(0_{X^{\ast }},0)\in \mathcal{F}^{\#}$.

$\left[ (\mathrm{iii})\Longrightarrow (\mathrm{i})\right] $ Assume that $(%
\mathfrak{A})$ holds. By $(H_{1})$ we have 
\begin{equation*}
\inf\limits_{x\in X}\sup\limits_{u\in U}(F_{u})^{\ast \ast }(x,0_{Y})\geq 0,
\end{equation*}%
that is $(0_{X^{\ast }},0)\in \limfunc{epi}q^{\ast \ast }$. By $(H_{1})$ we
have also $\limfunc{dom}q^{\ast }\neq \emptyset $. Consequently, $%
(0_{X^{\ast }},0)\in \overline{\limfunc{co}}\mathcal{F}^{\#}$ and, by (iii), 
$(0_{X^{\ast }},0)\in \mathcal{F}^{\#}$, that is $(\mathfrak{B})$.\hfill $%
\square \medskip $

\begin{corollary}
\label{cor:3-1} (\cite[Theorem 4.1]{VB-24}) Given $F\colon X\times
Y\longrightarrow \overline{\mathbb{R}}$, consider the following statements:

\textup{(i)} The $S_{F}$-procedure is valid.

\textup{(ii)} The $S_{F^{\ast \ast }}$-procedure is valid.

\textup{(iii)} $\left\{ (x^{\prime },s)\in X^{\ast }\times \mathbb{R}%
:\exists \mu \in Y^{\ast },\;F^{\ast }(x^{\prime },\mu )\leq s\right\} $ is
weak*-closed regarding $\{(0_{X^{\ast }},0)\}.$

We have $(\mathrm{i})\Longrightarrow (\mathrm{ii})\Longrightarrow (\mathrm{%
iii})$.

If, moreover, 
\begin{equation*}
(H_{2})\qquad \qquad \inf_{x\in X}F(x,0_{Y})=\inf_{x\in X}F^{\ast \ast
}(x,0_{Y})\neq +\infty
\end{equation*}%
holds, then $(\mathrm{i})\QTR{frametitle}{\Longleftrightarrow }(\mathrm{ii})%
\QTR{frametitle}{\Longleftrightarrow }(\mathrm{iii})$.
\end{corollary}

\textbf{Proof }Apply Theorem \ref{thm:3-1} to the case that $U=\{u\}$ is a
singleton and $F_{u}=F$. We have 
\begin{equation*}
\mathcal{F}^{\#}=\left\{ (x^{\prime },s)\in X^{\ast }\times \mathbb{R}%
:\exists \mu \in Y^{\ast },\;F^{\ast }(x^{\prime },\mu )\leq s\right\}
\end{equation*}%
that is the projection of $\limfunc{epi}F^{\ast }$ on $X^{\ast }\times 
\mathbb{R}$. Since the conjugate function $F^{\ast }$ is convex, the set $%
\mathcal{F}^{\#}$ is convex and the proof of Corollary \ref{cor:3-1} is
complete.\hfill $\square \medskip $

\begin{remark}
\label{rem:3-1} Condition $(H_{1})$ (resp. $(H_{2})$) is in particular
satisfied if $(F_{u})^{\ast \ast }(x,0_{Y})=F_{u}(x,0_{Y})$ for all $%
(u,x)\in U\times X$ (resp. $F^{\ast \ast }(x,0_{Y})=F(x,0_{Y})$ for all $%
x\in X$), and $\exists a\in X:\sup\limits_{u\in U}F_{u}(a,0_{Y})\neq +\infty 
$ (resp. $\exists a\in X:F(a,0_{Y})\neq +\infty $).
\end{remark}

\section{Robust generalized $S$-procedure with RHS function in $\Gamma (X)$}

In this section $X,Y$ are lcHtvs, $h\in \Gamma (X)$, and $F_{u}\colon
X\times Y\longrightarrow \overline{\mathbb{R}}$ for all $u\in U$. Let us
recall the statements introduced in Section \ref{sec:1}: 
\begin{align*}
& (\mathfrak{A}_{h})\quad \sup\limits_{u\in U}F_{u}(x,0_{Y})\geq h(x)\text{
for all }x\in X, \\
& (\mathfrak{B}_{h})\quad \text{ }\forall a^{\prime }\in \limfunc{dom}%
h^{\ast },\exists (\bar{u},\bar{\mu})\in U\times Y^{\ast }:(F_{\bar{u}%
})^{\ast }(a^{\prime },\bar{\mu})\leq h^{\ast }(a^{\prime }).
\end{align*}

\begin{lemma}
\label{lem:4-1} One has $(\mathfrak{B}_{h})\Longrightarrow (\mathfrak{A}%
_{h}) $.
\end{lemma}

\textbf{Proof }Assume $(\mathfrak{B}_{h}).$ For all $(x,a^{\prime })\in
X\times X^{\ast }$ there exists $\left( \bar{u},\bar{\mu}\right) \in U\times
Y^{\ast }$\ such that%
\begin{equation*}
\begin{array}{ll}
{\sup\limits_{u\in U}F_{u}(x,0_{Y})} & \geq F_{\bar{u}}(x,0_{Y})\geq (F_{\bar{u}%
})^{\ast \ast }(x,0_{Y}) \\ 
& \geq \langle a^{\prime },x\rangle +\langle \bar{\mu},0_{Y}\rangle -(F_{%
\overline{u}})^{\ast }(a^{\prime },\bar{\mu}) \\ 
& \geq \langle a^{\prime },x\rangle -h^{\ast }(a^{\prime }).%
\end{array}%
\end{equation*}%
Taking the supremum on all $a^{\prime }\in \limfunc{dom}h^{\ast }$ we obtain
that $(\mathfrak{A}_{h})$ holds.\hfill $\square \medskip $

The next result involves the graph of $h^{\ast }$, that is, $\limfunc{gph}%
h^{\ast }:=\left\{ {(a^{\prime },h^{\ast }(a^{\prime }))},\;a^{\prime }\in 
\limfunc{dom}h^{\ast }\right\} .$

\begin{theorem}
\label{thm:4-1} Assume that%
\begin{equation*}
(H_{3})\qquad 
\left\{  \begin{array}{ll}
\forall a^{\prime }\in \limfunc{dom}h^{\ast },\; & \inf\limits_{x\in
X}\left( \sup\limits_{u\in U}F_{u}(x,0_{Y})-\langle a^{\prime },x\rangle
\right) \\ 
& =\inf\limits_{x\in X}\left( \sup\limits_{u\in U}(F_{u})^{\ast \ast
}(x,0_{Y})-\langle a^{\prime },x\rangle \right) \neq +\infty%
\end{array}%
\right.
\end{equation*}%
and%
\begin{equation*}
 (H_{4})\qquad \mathcal{F}^{\#}\text{ is weak*-closed convex regarding }%
\limfunc{gph}h^{\ast }
\end{equation*}%
hold. Then, the robust $S_{\left\{ F_{u}\right\} _{u\in U}}^{h}$-procedure
is valid.
\end{theorem}

\textbf{Proof }One has to prove that $(\mathfrak{A}_{h})\Longrightarrow (%
\mathfrak{B}_{h})$. Let $a^{\prime }\in \limfunc{dom}h^{\ast }$. We have 
\begin{equation*}
q^{\ast \ast }(a^{\prime })=\sup_{x\in X}\left( \langle a^{\prime },x\rangle
-\sup_{u\in U}(F_{u})^{\ast \ast }(x,0_{Y})\right)
\end{equation*}%
and by $(H_{3})$, 
\begin{equation*}
q^{\ast \ast }(a^{\prime })=\sup_{x\in X}\left( \langle a^{\prime },x\rangle
-\sup_{u\in U}F_{u}(x,0_{Y})\right) .
\end{equation*}%
By $(\mathfrak{A}_{h})$ we obtain that 
\begin{equation*}
q^{\ast \ast }(a^{\prime })\leq \sup_{x\in X}\left( \langle a^{\prime
},x\rangle -h(x)\right) ,
\end{equation*}%
that is $(a^{\prime },h^{\ast }(a^{\prime }))\in \limfunc{epi}q^{\ast \ast }$%
. Since $\limfunc{dom}q^{\ast }\neq \emptyset $ (see $(H_{3})$), it ensures
that 
\begin{equation*}
(a^{\prime },h^{\ast }(a^{\prime }))\in \overline{\limfunc{co}}\mathcal{F}%
^{\#},
\end{equation*}%
and, by $(H_{4})$, 
\begin{equation*}
(a^{\prime },h^{\ast }(a^{\prime }))\in \mathcal{F}^{\#}\text{ for all }%
a^{\prime }\in \limfunc{dom}h^{\ast },
\end{equation*}%
that is $(\mathfrak{B}_{h})$.\hfill $\square \medskip $

We conclude with the determinist case when $U=\{u\}$, $F_{u}=F\colon X\times
Y\longrightarrow \overline{\mathbb{R}}$, $\mathcal{F}^{\#}=\left\{
(x^{\prime },s)\in X^{\ast }\times \mathbb{R}:\exists \mu \in Y^{\ast
},\;F^{\ast }(x^{\prime },\mu )\leq s\right\} $, and $h\in \Gamma (X)$. The
corresponding statements are 
\begin{align*}
& (\mathbb{A}_{h})\quad F(x,0_{Y})\geq h(x)\text{ for all }x\in X, \\
& (\mathbb{B}_{h})\quad \forall a^{\prime }\in \limfunc{dom}h^{\ast
},\exists \mu \in Y^{\ast }:F^{\ast }(a^{\prime },\mu )\leq h^{\ast
}(a^{\prime }).
\end{align*}

By Lemma \ref{lem:4-1} we have $(\mathbb{B}_{h})\Longrightarrow (\mathbb{A}%
_{h})$. The so-called $S_{F}^{h}$-procedure is said to be \textit{valid} if $%
(\mathbb{A}_{h})\QTR{frametitle}{\Longleftrightarrow }(\mathbb{B}_{h})$, or,
equivalently, $(\mathbb{A}_{h})\Longrightarrow (\mathbb{B}_{h})$.

\begin{corollary}
\label{cor:4-1} Assume that%
\begin{equation*}
(H_{5})\qquad \left\{ 
\begin{array}{ll}
\forall a^{\prime }\in \limfunc{dom}h^{\ast },\; & \inf_{x\in X}\left(
F(x,0_{Y})-\langle a^{\prime },x\rangle \right) \\ 
& =\inf\limits_{x\in X}\left( F^{\ast \ast }(x,0_{Y})-\langle a^{\prime
},x\rangle \right) \neq +\infty%
\end{array}%
\right.
\end{equation*}%
and%
\begin{equation*}
(H_{6})\qquad 
\begin{cases}
\text{The \ set }\left\{ (x^{\prime },s)\in X^{\ast }\times \mathbb{R}%
:\exists \mu \in Y^{\ast },\;F^{\ast }(x^{\prime },\mu )\leq s\right\} \text{
} \\ 
\text{is weak*-closed regarding }\limfunc{gph}h^{\ast }%
\end{cases}%
\end{equation*}%
hold. Then, the $S_{F}^{h}$-procedure is valid.
\end{corollary}

\textbf{Proof }Apply Theorem \ref{thm:4-1}, noticing that $\mathcal{F}^{\#}$
is convex.\hfill $\square \medskip $

\begin{remark}
\label{rem:4-1} Condition $(H_{5})$ is in particular satisfied if $%
F(x,0_{Y})=F^{\ast \ast }(x,0_{Y})$ for all $x\in X$ and there exists $a\in
X $ such that $F(a,0_{Y})\neq +\infty $.
\end{remark}

{\bf Acknowledgements}  This work is partly supported by the project “Generalized Farkas  lemma for a  family of adjustable systems with uncertainty and applications",  Vietnam National University-Ho Chi Minh city,
Vietnam.

\end{document}